\documentclass[12pt,reqno]{amsart}

\usepackage{pb-diagram}
\usepackage{mathrsfs}

\newtheorem{theorem}{Theorem}[section]

\theoremstyle{definition}
\newtheorem{definition}[theorem]{Definition}

\theoremstyle{remark}

\numberwithin{equation}{section}

\author{Gizem Karaali}
\address{Department of Mathematics, Pomona College, Claremont CA
91711}
\email{gizem.karaali@pomona.edu}

\subjclass[2000]{Primary 17B62, 17B20} 

\date{August 18, 2007.}

\keywords{Hopf algebra, quantum group, quasi-Hopf algebra, weak Hopf algebra,
Hopf algebroid, Hopf monad, hopfish algebra}

\theoremstyle{definition}
\newtheorem*{acknow}{Acknowledgments}

\newcommand{\baseRing}[1]{\ensuremath{\mathbb{#1}}}
\newcommand{\C}{\baseRing{C}}

\newcommand{\Z}{\baseRing{Z}}

\newcommand{\kk}{\baseRing{K}}
\newcommand{\g}{\mathfrak{g}}

\newcommand{\h}{\mathfrak{h}}

\renewcommand{\phi}{\varphi}

\allowdisplaybreaks

\begin{document}

\title{On Hopf algebras and their Generalizations}

\begin{abstract}
We survey Hopf algebras and their generalizations. In particular, we compare and
contrast three well-studied generalizations (quasi-Hopf algebras, weak Hopf
algebras, and Hopf algebroids), and two newer ones (Hopf monads and hopfish algebras). Each of these notions was originally introduced for a specific
purpose within a particular context; our discussion favors applicability to the
theory of dynamical quantum groups. Throughout the note, we provide several
definitions and examples in order to make this exposition accessible to readers
with differing backgrounds.
\end{abstract}

\maketitle

\section{Introduction - Goals and Motivation}

The purpose of this note is to provide a historical and comparative study of the
several notions of generalized Hopf algebras that have been 
introduced and studied in recent years. We start with a brief discussion of Hopf
algebras. We then consider the three main contenders for the \emph{correct}
notion of a generalized Hopf algebra: quasi-Hopf algebras, weak Hopf algebras,
and Hopf algebroids. We then discuss two newer notions that generalize Hopf algebras: Hopf monads and hopfish algebras. As these are newer concepts, our study of them is necessarily rather cursory.

In our pursuit we trace the steps of many researchers, both mathematically and
philosophically. Hence, this is an expository note, with no claim of
introducing new mathematics. Its sole 
purpose is to collect together several pieces of interesting mathematics and
present them in one comprehensive historical narrative, so 
as to provide a broader perspective on the current status of research involving
various generalizations of Hopf algebras. 

The technical details and basic examples provided are intended to make this note
accessible to a beginner in the theory of Hopf algebras, while we hope that
those with more experience will still enjoy reading the discussions. 

\begin{acknow}
The author would like to thank Marcelo Aguiar, Gabriella B\"{o}hm, Ignacio Lopez Franco, Ludmil Hadjiivanov, Mitja Mastnak, Susan Montgomery, Frederic Patras, Fernando Souza, Ivan Todorov, Alexis Virelizier, Alan Weinstein and Milen Yakimov for helpful comments and useful questions during various stages of the work that led to this paper.
\end{acknow}

\section{Hopf algebras}
\label{HopfAlgebraDefinitionsSection}

Naturally, we begin with the definition of a Hopf algebra:

\begin{definition}
\label{Hopfdef}
A \textbf{Hopf algebra} over a commutative ring $\kk$ is a $\kk$-module $H$ such
that:
\begin{enumerate}
\item $H$ is an associative unital $\kk$-algebra (with
${m : H \otimes H \rightarrow H}$ as the \emph{multiplication} and 
${u : \kk \rightarrow H}$ as the \emph{unit}) and a coassociative counital
$\kk$-coalgebra (with
${\Delta : H \rightarrow H \otimes H}$ as the \emph{comultiplication} and 
${\epsilon : H \rightarrow \kk}$ as the \emph{counit});
\item The comultiplication $\Delta$ and the counit ${\epsilon}$ are both
algebra homomorphisms;
\item The multiplication ${m}$ and the unit ${u}$ are both coalgebra
homomorphisms;
\item There is a bijective $\kk$-module map ${S : H \rightarrow H}$,
called the \emph{antipode}, such that for all elements ${h \in H}$:
\begin{equation}
\label{antipodedef}
S(h_{(1)})h_{(2)} = (u \circ \epsilon) (h) = h_{(1)}S(h_{(2)}). 
\end{equation}
\end{enumerate}
\end{definition}

We first note here that in the defining formula
for the antipode, we are making use of the famous \emph{Sweedler notation}\footnote{In the introduction to \cite{Swe69}, Sweedler remarks that this notation was developed during years of joint work with Robert Heyneman. However, it is much more common to see references to Sweedler notation than to Heyneman-Sweedler notation. We follow the tradition here.},
which is quite well-established in the Hopf algebra literature. 
In short the Sweedler notation is a generalization of the Einstein notation (in that it intrinsically demands a summation). More specifically for any $h \in H$, we
write:
\[ \Delta(h) = h_{(1)} \otimes h_{(2)}. \]
\noindent
This presentation itself is purely symbolic; the terms $h_{(1)}$ and $h_{(2)}$
do not stand for particular elements of $H$. The comultiplication $\Delta$
takes values in $H \otimes H$, and so we know that:
\[ \Delta(h) = (h_{1,1} \otimes h_{1,2}) + (h_{2,1} \otimes h_{2,2}) +
(h_{3,1}\otimes h_{3,2}) + \cdots + (h_{N,1}\otimes h_{N,2}) \]
\noindent
for some elements $h_{i,j}$ of $H$ and some integer $N$. The Sweedler notation
is just a way to separate the $h_{i,1}$ from the $h_{j,2}$. In other words,
one can say that the notation $h_{(1)}$ stands for the generic $h_{i,1}$ and
the notation $h_{(2)}$ stands for the generic $h_{j,2}$. However the summation
is inherent in the notation: Whenever one sees a term like $S(h_{(1)})h_{(2)}$
(the left hand side of \eqref{antipodedef}), one has to realize that this
stands for a sum of the form: 
\[S(h_{1,1}) \otimes h_{1,2} + S(h_{2,1})\otimes h_{2,2} + S(h_{3,1})\otimes
h_{3,2} + \cdots + S(h_{N,1})\otimes h_{N,2}. \]

A second remark that should be made here involves commutative diagrams. In most
texts introducing the notion of a Hopf algebra, one will come across
commutative diagrams similar to the following:
\[ \begin{diagram}
     \node{H\otimes H} \arrow{e,t}{S \otimes id}  
         \node{H \otimes H} \arrow{s,r}{m} \\
     \node{H} \arrow{n,l}{\Delta} \arrow{e,b}{u \circ \epsilon} \node{H}
   \end{diagram}
\qquad \qquad
\begin{diagram}
    \node{H\otimes H} \arrow{e,t}{id \otimes S}  
         \node{H \otimes H} \arrow{s,r}{m} \\
     \node{H} \arrow{n,l}{\Delta} \arrow{e,b}{u \circ \epsilon} \node{H}
   \end{diagram}   \]
In fact these two diagrams are equivalent to the single-line equation
\eqref{antipodedef}, but for those comfortable with commutative diagrams, they
provide a much more visual way to understand it. Incidentally an even more
compact way to describe $S$ would be as the convolution inverse of the
identity homomorphisn on $H$. 

For the sake of completeness, we include here the diagrams representing the
properties of the counit (the so-called \emph{counit axiom}):
\[ \begin{diagram}
\node{ }\node{H} \arrow{sw,l}{\Delta} \arrow{se,l}{\Delta} \arrow[2]{s,r}{id_H}
\node{ } \\
\node{H \otimes H} \arrow{se,r}{\epsilon \otimes 1} \node{ } \node{H \otimes H}
\arrow{sw,r}{1 \otimes \epsilon}\\
\node[2]{H}
\end{diagram}\]
\noindent
(where we use $H \cong H \otimes \kk \cong \kk \otimes H$), and the
comultiplication (the coassociativity): 
\[ \begin{diagram}
     \node{H} \arrow{e,t}{\Delta} \arrow{s,l}{\Delta} 
         \node{H \otimes H} \arrow{s,r}{\Delta \otimes 1}  \\
     \node{H \otimes H} \arrow{e,b}{1 \otimes \Delta} \node{H \otimes H \otimes
H}
   \end{diagram}\]   

We also should point out that the second and third conditions of Definition
\ref{Hopfdef} are equivalent to one another and, together with the first
condition, define what is called a \emph{$\kk$-bialgebra}:
\begin{definition}
\label{bialgebradef}
A \textbf{bialgebra} over a commutative ring $\kk$ is a $\kk$-module $B$ such
that:
\begin{enumerate}
\item $B$ is both an associative unital $\kk$-algebra (with
${m : B \otimes B \rightarrow B}$ as the \emph{multiplication} and 
${u : \kk \rightarrow B}$ as the \emph{unit}) and a coassociative counital
$\kk$-coalgebra (with
${\Delta : B \rightarrow B \otimes B}$ as the \emph{comultiplication} and 
${\epsilon : B \rightarrow \kk}$ as the \emph{counit});
\item The comultiplication $\Delta$ and the counit ${\epsilon}$ are both
algebra homomorphisms.
\end{enumerate}
\end{definition}
\noindent
Then we can define a Hopf algebra to be a bialgebra with an antipode, in other
words, a bialgebra $B$ with an antihomomorphism $S : B \rightarrow B$
satisfying Equation
\eqref{antipodedef}.
   
Hopf algebras were first introduced in the early 1940s by Hopf in \cite{Ho41},
where he was working on homology rings of certain compact manifolds. (For a
modern exposition of his results using the language of Hopf algebras, we refer
the reader to \cite{McC}). Later on, more examples were discovered and basic
references started appearing, see for instance the early classics like
Sweedler's \cite{Swe69} and Abe's \cite{Ab80}. For a more modern approach
emphasizing actions of Hopf algebras, see Montgomery's \cite{Mo93}.

For a simple example of Hopf algebras, consider the group algebra $\kk[G]$ of a
group $G$ over a field $\kk$. Here we define the coalgebra structure by
defining, on the elements $g$ of $G$, the comultiplication and the counit by
setting $\Delta(g) = g \otimes g$ and $\epsilon(g) = 1$. (More generally, for
an arbitrary coalgebra $C$, the elements $g$ of $C$ satisfying $\Delta(g) = g
\otimes g$ and $\epsilon(g) = 1$ are called \emph{group-like}). The antipode
is defined on the group elements by $S(g)= g^{-1}$.

Dually, we can look at the algebra $F_{\kk}(G)$ of $\kk$-valued functions on a
group
$G$.\footnote{We are intentionally being vague, in order to include several
classes of algebras here; if, for instance, $G$ is an (affine) algebraic
group, we could be thinking of the coordinate ring of $G$; if $G$ is a Poisson
group, we could be looking at the Poisson algebra $C^{\infty}(G)$.} 
There is a natural (commutative) multiplication $m$ on this algebra given by: 
\[ m(f_1,f_2) (g) = f_1(g) f_2(g) \textmd{ for all } g \in G. \]
\noindent
The unit just maps elements of the field $\kk$ to the associated constant
functions. To define a comultiplication, one needs to use the natural
embedding of $F_{\kk}(G) \otimes F_{\kk}(G)$ into $F_{\kk}(G \times G)$. Then
the group
multiplication $m_G : G \times G \rightarrow G$ on $G$ induces a
comultiplication as follows:
\[ \Delta(f) (g_1,g_2) = f(m_G(g_1,g_2)). \]
\noindent
The counit is the map taking $f \in F_{\kk}(G)$ to $f(e_G) \in {\kk}$ where
$e_G$ is
the identity element of the group $G$. Finally the antipode is the map $S$
defined by $S(f)(g) = f(g^{-1})$. 

Function algebras give us a large collection of examples of \emph{commutative}
Hopf algebras, (which are just Hopf algebras with commutative multiplications).
In fact, one can show that any finitely generated commutative Hopf algebra over
a field $\kk$ of characteristic zero is the function algebra of an affine
algebraic group $G$ over $\kk$. 
We will come back to this in Section \ref{HopfMonadSection}.

Another family of examples is given by universal enveloping algebras $U(\g)$ of
Lie algebras $\g$. In this case, the coalgebra structure is defined on the
elements $x$ of $\g$ by:
\[ \Delta(x) = x \otimes 1 + 1 \otimes x; \qquad \qquad \epsilon(x) = 0. \]
\noindent
(More generally, for an arbitrary coalgebra $C$, elements $x$ of $C$ satisfying 
$\Delta(x) = x \otimes 1 + 1 \otimes x$ and $\epsilon(x) = 0$ are called
\emph{primitive}). The antipode $S$ is defined on the Lie algebra elements as
$S(x) = -x$.

This last class of examples provides us with a large collection of
\emph{cocommutative} Hopf algebras. Cocommutativity is the property of the
comultiplication $\Delta$ described by the following commutative diagram:
\[  \begin{diagram}
    \node{H\otimes H} \node{H \otimes H} \arrow{w,t}{T}  \\
     \node{H} \arrow{n,l}{\Delta} \node{H} \arrow{n,r}{\Delta} \arrow{w,b}{id} 
   \end{diagram} \]
\noindent
which is exactly the diagram one obtains by switching the directions of the
arrows and replacing $m$ by $\Delta$ in the diagram describing commutativity
of a product $m$:
\[  \begin{diagram}
    \node{H\otimes H} \arrow{s,l}{m} \arrow{e,t}{T} \node{H \otimes H}
\arrow{s,r}{m} \\
     \node{H}  \arrow{e,b}{id}  \node{H} 
   \end{diagram} \]
\noindent
In both of these diagrams $T: H \otimes H \rightarrow H \otimes H$ is the
usual twist map: $T(a \otimes b) = b \otimes a$. 

After seeing the result above about commutative Hopf algebras, one may expect
an analogous result regarding cocommutative Hopf algebras. In fact one can show
that 
a cocommutative Hopf algebra generated by its primitive elements is indeed the
universal enveloping algebra of a Lie algebra. 
This result will be mentioned in Section \ref{HopfMonadSection} as well.

Ever since their introduction, Hopf algebras have been studied by many
mathematicians. In the early 1970s, Hochschild, while developing the theory of
algebraic groups, translated much of representation theory into the language of
Hopf algebras. (See for instance \cite{Ho81} for a classic written with this
perspective; \cite{FeRi05} is a more modern text using the same approach). In
some sense, we can say that Hopf algebras provided
mathematicians the ultimate framework to do representation theory.

The development of the theory of quantum groups brought a fresh revival of
interest in Hopf algebras.\footnote{In the following
paragraphs, we will provide only a brief sketch of the ideas of quantum group
theory and refer the more interested reader to one of many textbooks and
monographs in the subject (eg. \cite{BrGo, CP, EtScbook, HoKa, Ka95, LaRa,
LusBook, Maj95, SS93}).}
Since the early 1980s, many mathematicians have
been working on structures which today we loosely call \emph{quantum groups}.
The earliest examples of quantum groups were particular deformations of
universal enveloping algebras of simple Lie algebras and function algebras of
simple algebraic groups. In mid 1980s, Drinfeld showed that the correct
framework to use when studying quantum groups was that of Hopf algebras; see,
for instance, his ICM address \cite{D1}. Thus the known world of Hopf algebras
was significantly expanded to include all these new examples, and results from
quantum groups began to add more spice and flavor to the classical theory.

Broadly speaking, a quantum group is a special type of noncommutative noncocommutative Hopf algebra. One obtains such Hopf algebras by deforming the multiplication or the
comultiplication of a commutative or a cocommutative
Hopf algebra. Since the end result of such a deformation is not commutative, it cannot properly be associated to a group and be its function algebra. Similarly, since it is not cocommutative, it is by no means the enveloping algebra of a Lie algebra. However, we can still view it as if it
were the function algebra or the enveloping algebra of some phantom group or Lie
algebra, and thus we have the term ``\textit{quantum group}".\footnote{This may remind some readers the philosophy of noncommutative
geometry a la Connes. A similar approach to quantum groups would involve
viewing them as symmetry objects of some quantum space; see \cite{Man88}.}

For the sake of completeness we briefly describe the most well-known quantum
group here: \textit{Quantum $sl_2(\C)$}, denoted by $U_h(sl_2(\C))$. The
notation is more illuminating than the name; this is a particular (Hopf
algebra) deformation of the universal enveloping algebra $U(sl_2(\C))$ of
$sl_2(\C)$. 

As an algebra, it is generated by $E, F, H$ subject to the following relations:
\[  [H,E] = 2E, \qquad [H, F] = -2F, \qquad [E,F] = \frac{q^H - q^{-H}}{q -
q^{-1}}\]
\noindent
where $h$ is viewed as a formal parameter, $q = e^h = \sum_{n=0}^{\infty}
\tfrac{h^n}{n!} \in \C[[h]]$, and $q^H, q^{-H}$ should be interpreted in a
similar manner. The coalgebra structure is defined on the generators by:
\begin{align*}
\Delta(E) &= E \otimes q^H + 1 \otimes E; \quad &\epsilon(E) = 0; \\
\Delta(F) &= F \otimes 1 + q^{-H} \otimes F; \quad &\epsilon(F) = 0; \\
\Delta(H) &= H \otimes 1 + 1 \otimes H; \quad &\epsilon(H) = 0.
\end{align*}
\noindent
Finally we define the antipode on the generators as:
\[ S(E) = -Eq^{-H}; \qquad S(F) = -q^HF; \qquad S(H) = -H.\]
Extended linearly, these give us a Hopf algebra structure, the famous quantum
$sl_2$. 

Currently, there is no consensus on what the precise definition of a quantum
group should be. Mainly,
there are several well-known examples and an exponentially growing literature
investigating their properties. For the purposes of this note, it will suffice
to 
identify the term quantum group with quasitriangular or co-quasitriangular 
Hopf algebras, i.e. noncommutative noncocommutative Hopf algebras 
associated to solutions of a certain equation, the quantum Yang-Baxter 
Equation, which we will describe in the next section.

\section{Long Interlude - Why generalize?}
\label{InterludeSection1}

We believe that the previous section provides a sufficient overview
of Hopf algebras, and prepares the reader for the discussion on the
generalizations of this classical theory. For the reader who wishes to
learn more about Hopf algebras, we recommend the references already mentioned
above, as well as any textbook on quantum groups.\footnote{Among the many books on
quantum groups which also have a detailed exposition of Hopf algebras, we list
here only a few: \cite{CP}, \cite{EtScbook}, \cite{Ka95}, \cite{Maj95},
\cite{SS93}.
A quick introduction to Hopf algebras and their applications may be found in \cite{Haz91}. Other references about Hopf algebras which especially 
focus on combinatorial applications are \cite{NiSw, AgMa}.}

Before we move on to the comparative study and technical details of the
various generalizations, we would like to philosophize a bit about why we are
interested in \emph{any} generalization. To those who view mathematics as an
intellectual pursuit merely interested in pure abstractions, the answer will be
clear: \emph{Why not?} However, for those who may need more motivation, we will
provide one, which comes from the theory of quantum groups.

The study of quantum groups goes back to the well-known \emph{quantum
Yang-Baxter equation}, which in its simplest form is as follows:
\begin{equation}
\label{QYBE}
R^{12}R^{13}R^{23} = R^{23}R^{13}R^{12}.
\end{equation}
\noindent
Here we can view $R$ as a map ${R : G \times G \rightarrow G \times G}$ for
some factorizable Poisson-Lie group (i.e. a Lie group with a Poisson bracket
on its function algebra). Then:
\[ R^{12}(g,h,k) = (R(g,h),k),\]
\noindent
and $R^{13}$ and $R^{23}$ can be defined likewise. In the realm of Hopf
algebras, the \emph{quantum Yang-Baxter equation} (henceforth referred to
as the QYBE) is the same equation \eqref{QYBE}; however, this time its 
solutions, the so-called \emph{quantum $R$-matrices}, are linear maps of 
$\kk$-modules: ${R \in H \otimes H}$. For surveys on the quantum Yang-Baxter
equation, displaying the several connections to the physics literature, see 
\cite{Ji1, Ji2}.

The QYBE and its solutions give rise to quantum groups: in the geometric
(Poisson) picture, they give a (Hopf algebra) deformation of the function
algebra of the relevant group $G$; in the Lie algebra picture, they give a
(Hopf algebra) deformation of the associated universal enveloping algebra. 
In the case of quantum $sl_2$ as presented above, the relevant quantum
$R$-matrix is: 

\[ R = q^{\tfrac{1}{2}H\otimes H} \sum_{n \ge 0} q^{\tfrac{n(n-1)}{2}}
\frac{[q-q^{-1}]^n}{[n]_q!} E^n \otimes F^n.\]
\noindent
Here we are using the $q$-notation common in the literature:
\[ [n]_q = \frac{q^n - q^{-n}}{q-q^{-1}} = q^{n-1} + q^{n-2} + \cdots + q^{-n+1}
\textmd{ and }
[n]_q! = \Pi_{k=1}^n [k]_q \]

In 1984 \cite{GN}, a modified version of the QYBE appeared in the mathematical physics
literature.
This new equation, later named the
\emph{quantum  dynamical Yang-Baxter equation}, (henceforth labeled QDYBE),  
lives in ${V \otimes V \otimes V}$ for a semisimple module $V$ of an abelian Lie 
algebra $\h$: 
\begin{equation}
\label{QDYBE}
R^{12}(\lambda-h^{(3)})R^{13}(\lambda)R^{23}(\lambda - h^{(1)}) =
R^{23}(\lambda)R^{13}(\lambda-h^{(2)})R^{12}(\lambda).
\end{equation}
\noindent
Here, the term $h^{(i)}$ is to be substituted by $\mu_i$ if $\mu_i$ is the weight of
the $i$th tensor component, $i= 1,2,3$. 
A solution ${R : \h^* \rightarrow End_{\h}(V \otimes V)}$ to the QDYBE is
called a \emph{quantum dynamical $R$-matrix}. For an introduction to the
dynamical Yang-Baxter equations along with some of their solutions, we refer
the reader to \cite{Fe1}. A more geometric exposition can be found in
\cite{EV2}.

Studying the QDYBE and its solutions, we get into the realm of \emph{dynamical
quantum groups}. The first examples of dynamical quantum groups that appeared
in the literature are Felder's \textit{elliptic quantum groups}, which were
introduced in \cite{Fe1, Fe2}.  A standard example that is studied in much
detail is $E_{\tau,\eta}(sl_2)$, which is an algebra over $\C$ generated by two
kinds of generators: meromorphic functions of a single variable $h$, and matrix
elements of a matrix $L(\lambda, w) \in End(\C^2)$. The two subscripts $\tau,
\eta$ are nonzero complex numbers with $Im(\tau) > 0$. 

Without going into much detail, we note that the $4 \times 4$ matrix solution to
the QDYBE associated with $E_{\tau,\eta}(sl_2)$ is in the following form:
\begin{eqnarray*}
R(\lambda,w,\tau,\eta) &=& E_{1,1} \otimes E_{1,1} + E_{2,2} \otimes E_{2,2}  \\
&+& \alpha(\lambda,w,\tau,\eta)E_{1,1} \otimes E_{2,2} +
\beta(\lambda,w,\tau,\eta)E_{1,2} \otimes E_{2,1} \\
&+& \gamma(\lambda,w,\tau,\eta)E_{2,1} \otimes E_{1,2} +
\delta(\lambda,w,\tau,\eta)E_{2,2} \otimes E_{1,1}
\end{eqnarray*}
\noindent
Here the functions $\alpha,\beta,\gamma,\delta$ are defined in terms of the
theta function:
\[ \theta(z,\tau) = - \Sigma_{j=-\infty}^{\infty} \exp(\pi i (j +
\tfrac{1}{2})^2\tau + 2\pi i(j+\tfrac{1}{2})(z + \tfrac{1}{2}))\]
\noindent
by:
\begin{eqnarray*}
\alpha(\lambda,w,\tau,\eta) =
\frac{\theta(w)\theta(\lambda+2\eta)}{\theta(w-2\eta)\theta(\lambda)} &,& 
\beta(\lambda,w,\tau,\eta) =
\frac{\theta(-w-\lambda)\theta(2\eta)}{\theta(w-2\eta)\theta(\lambda)} \\
\gamma(\lambda,w,\tau,\eta) =
\frac{\theta(w-\lambda)\theta(2\eta)}{\theta(w-2\eta)\theta(\lambda)} &,&
\delta(\lambda,w,\tau,\eta) =
\frac{\theta(w)\theta(\lambda-2\eta)}{\theta(w-2\eta)\theta(\lambda)}
\end{eqnarray*}
\noindent
For more on $E_{\tau, \eta}(sl_2)$, including an explicit presentation in terms
of generators and relations, we refer the reader to \cite{FeVa1}. 

The objects, like Felder's elliptic quantum groups, associated to solutions of
the QDYBE turn out to be quite Hopf-like in many respects, but they are not all
necessarily Hopf algebras. Or in other words, we can say that, just as in the
case of quantum groups, some Hopf-like structures come into play when one
studies the solutions of the QDYBE. And herein lies our motivation for the
purpose of this particular note. A need for the \emph{correct Hopf-like object}
that will provide the framework for the theory of dynamical quantum groups
leads us to the study of various generalizations of Hopf algebras. 

\section{Candidate 1 : Quasi-Hopf Algebras}

Our first candidate in our search for the correct generalization of Hopf 
algebras is the quasi-Hopf algebra, first introduced by Drinfeld in \cite{D4,
D5}, and used 
to give a natural proof of Kohno's theorem relating the monodromy of the
Knizhnik-Zamolodchikov equations to a representation of the braid group arising
from a quantum group. 
In order to define quasi-Hopf algebras, we first need to introduce the notion of
a \emph{quasi-bialgebra}:

\begin{definition}
\label{quasibialgebradef}
A \textbf{quasi-bialgebra} $B$ over a commutative ring $\kk$ is a unital
associative $\kk$-algebra equipped with two algebra homomorphisms ${\epsilon :
B \rightarrow \kk}$ (the \emph{counit}) and 
${\Delta : B \rightarrow B \otimes B}$  (the \emph{comultiplication}) together
with an invertible element $\Phi$ of $B \otimes B \otimes B$. Furthermore we
require the following to hold:
\begin{enumerate}
\item $\Phi$ satisfies the \emph{pentagon relation}:
\[ (\Delta \otimes 1 \otimes 1)\Phi \cdot (1 \otimes 1 \otimes \Delta) \Phi =
(\Phi \otimes 1) \cdot (1 \otimes \Delta \otimes 1)\Phi \cdot (1 \otimes
\Phi)\]
\item $\Delta$ is \emph{quasi-coassociative}:
\[ (1 \otimes \Delta) \Delta(b) = \Phi^{-1} (\Delta \otimes 1) \Delta(b) \Phi
\textmd{ for all } b \in B \] 
\item $\Phi$ is compatible with the counit $\epsilon$ in the following sense:
\[ (1 \otimes \epsilon \otimes 1) \Phi = 1\]
\item $\epsilon$ satisfies the counit axiom.\footnote{We presented the
commutative diagram version of the counit axiom earlier, in closed form it
reads: 
\[ (\epsilon \otimes 1) \Delta = 1 = (1 \otimes \epsilon) \Delta. \]}
\end{enumerate}
\end{definition}

From (1), (3) and (4), we get: 
\[ (\epsilon \otimes 1 \otimes 1) \Phi = (1 \otimes \epsilon \otimes 1) \Phi =
(1 \otimes 1 \otimes \epsilon) \Phi = 1.\]

$\Phi$ of this definition is alternatively called an \emph{associator} or a
\emph{co-associator}. We will call it the associator, merely to ease our typing
efforts. 
However we should note that both versions have merit. The term
``\textit{co-associator}" makes sense because $\Phi$ gives information on how
far $\Delta$ is from coassociativity. The term ``\textit{associator}" makes
sense because using $\Phi$, one can define associativity isomorphisms  in the
representation category of $B$ which makes it into a (non-strict) monoidal
category. The pentagon relation (1) satisfied by $\Phi$ in the above definition
translates automatically to the pentagon relations of the associativity
isomorphisms in the representation category of $B$.

Comparing Definitions \ref{bialgebradef} and \ref{quasibialgebradef}, it is
clear that this generalization has cost us only the coassociativity of the
comultiplication. However this price is perfectly acceptable to many
researchers, as quasibialgebras and their Hopf relatives fit very nicely into
several theories. We will discuss these shortly.

Here is the definition of the Hopf version:\footnote{We note here that our
definitions  in this section are not quite as general as the original ones in
\cite{D5}. Drinfeld's original definitions involved two invertible elements $l$
and $r$, or equivalently a left unit constraint $l$ and a right unit constraint
$r$, satisfying the so-called \emph{Triangle Axiom}. However, Drinfeld showed,
also in \cite{D5}, that he could always reduce his quasi-Hopf algebras into
quasi-bialgebras where $r=l=1$. Therefore we will not be worried much about our
more restrictive definitions.}

\begin{definition}
\label{quasiHopfdef}
A \textbf{quasi-Hopf algebra} $H$ over a commutative ring $\kk$ is a
quasi-bialgebra over $\kk$ equipped with an (algebra and coalgebra)
antihomomorphism $S : H \rightarrow H$ and two 
canonical elements $\alpha, \beta \in H$ such that:
\begin{enumerate}
\item $m \cdot (S \otimes \alpha) \Delta(b) = \epsilon(b) \alpha$ for all $b \in
B$;
\item $m \cdot (1 \otimes \beta S) \Delta(b) = \epsilon(b)\beta$ for all $b \in
B$;
\item $m (m \otimes 1) \cdot (S \otimes \alpha \otimes \beta S) \Phi = 1$; 
\item $m (m \otimes 1) \cdot (1 \otimes \beta S \otimes \alpha) \Phi^{-1} = 1$.
\end{enumerate}
\end{definition}

We note that (1) and (2) are the modified versions of the antipode axiom
\eqref{antipodedef}.
If $\alpha = \beta = 1$, and $\Phi = 1 \otimes 1 \otimes 1$, $H$ clearly becomes
a Hopf algebra.

Drinfeld introduced quasi-Hopf algebras in his work relating the
Knizhnik-Zamolodchikov equations to the theory of quantum groups. However,
these structures proved to be a lot more than mere tools in one proof..
The large interest they invoked in the mathematicians in the field can easily be deduced, for example, from J. Stasheff's enthusiastic MathSciNet review $MR1091757$ of Drinfeld's \cite{D5}.

One general reason for this excitement can be explained, in category-theoretic
terms, as follows: Quasi-Hopf algebras have the right kind of representation categories. In particular their representation categories are (not necessarily strict) rigid monoidal categories.\footnote{Here we digress briefly and venture into category theory. In Section \ref{HopfMonadSection}, we will give more details. The
classic reference for the terms we use is \cite{MacBook}. \cite{Ka95} provides
ample background and all the details within the context of Hopf algebras and
quantum groups. Note that in  \cite{Ka95}, our monoidal categories are called \emph{tensor categories}. These are exactly the ``\emph{categories with multiplication}" of \cite{Ben63, Mac63}. Also see the historical notes at the end of Chapter VII in \cite{MacBook}.}

Recall that a monoidal category is \textit{strict} if its associativity isomorphism is canonical. It can be shown that any (nice enough) strict monoidal category is the
representation category of a bialgebra, and any (nice enough) rigid strict monoidal category is the representation category of a Hopf algebra.\footnote{Obviously the phrase ``\emph{nice enough}" needs to be explained for this vague sentence to become a precise (and correct) mathematical statement. For instance, for a Tannaka-Krein type construction (cf. Section \ref{HopfMonadSection}), we need to construct a fiber functor from our category into the category of vector spaces. We will not attempt this here.} 
A well-known theorem of MacLane asserts that every non-strict monoidal category is in fact tensor equivalent to a strict one. So it is natural to look at algebraic objects whose representation categories are non-strict monoidal because these have in essence the same type of representation theory as that of bialgebras and Hopf algebras. 

It turns out that quasi-Hopf algebras are just the right structures in this
perspective! More specifically, quasi-bialgebras are precisely those algebras,
equipped with counit and comultiplication, whose representation categories are
monoidal. In fact in \cite{Ka95}, quasi-bialgebras are defined to be those structures which have monoidal representation categories. (This is immediately followed by a proof of the equivalence of this definition to the definition of Drinfeld in \cite{D5}). Likewise quasi-Hopf algebras are structures with rigid monoidal representation categories.

When the first examples of dynamical quantum groups began to appear, people
noticed that they were not Hopf algebras, but still looked very much like them.
Therefore several researchers focused on quasi-Hopf algebras, as structures
already under serious inspection for the reasons mentioned earlier, with the
expectation that these could possibly provide the right Hopf-like theory to
describe dynamical quantum groups. Indeed, it turns out that Felder's elliptic
quantum groups naturally fit into the framework of quasi-Hopf algebras,
\cite{EnFe1}.

Besides Felder's elliptic quantum groups, the fundamental example of a
nontrivial quasi-Hopf algebra (i.e. one that is not a Hopf algebra) is the one
constructed by Drinfeld. This construction involves the monodromy of the
Knizhnik-Zamolodchikov equations and is beyond the scope of this note. We refer
the interested reader to the original work of Drinfeld in \cite{D5} and the
more pedagogical exposition in \cite{EtScbook}. 

The representation theory of quasi-Hopf algebras is quite exciting.
Category-theoretic arguments ensure that quasi-objects will generally be only a
twist away from their non-quasi counterparts.\footnote{In this context, several
authors prefer to use terms like ``\textit{gauge transformation}" or
``\textit{skrooching}" in place of the term ``\textit{twist}."}
More precisely, MacLane's theorem about the (tensor-)equivalence of non-strict
monoidal categories to strict ones implies that the representation theory of a
quasi-object is quite similar to that of a non-quasi one. 

Tensor equivalence in the realm of monoidal categories translates into the
language of equivalence of quasi-bialgebras and quasi-Hopf algebras in terms of
twists in the algebraic realm. In many cases (when a certain cohomology
vanishes), the twist will involve a trivial isomorphism on the level of the
algebra. In other words, the quasi-structure will be ``\textit{twistable}" into
a non-quasi one on the same underlying algebra. In other cases we will need to
modify the underlying set. Once again the reader is referred to
\cite{EtScbook}, for a discussion of the use of twists of quasi-bialgebras and
quasi-Hopf algebras.

Overall, quasi-Hopf algebras are interesting structures, with many researchers
still investigating their various applications. Their representation theory is
also very appealing. Moreover as we mentioned earlier, they can be used to
describe the elliptic quantum groups of Felder.  Nonetheless, we continue our
search and look for other alternatives. This is mainly due to the general fact
that algebraists are not very much delighted by a non-associative or a
non-coassociative structure. The relaxation of these properties makes
description of actions and coactions a lot more cumbersome, and this is
typically not very desirable. Another issue is the notion of duality for
quasi-objects. Since the definition of a quasi-bialgebra is not self-dual (the
underlying set is an associative algebra but the coalgebra structure is not
coassociative), the 
natural object that should be the dual of a quasi-bialgebra is not a
quasi-bialgebra; and similarly the natural object which should be the dual of a
quasi-Hopf algebra is not a quasi-Hopf algebra.
One needs to define separately a dual quasi-bialgebra and hence a dual
quasi-Hopf algebra; see for instance the preliminaries of \cite{Maj92}.  As a
consequence, we do not have a natural notion of a comodule or a Hopf module
over a quasi-object and any other mathematical construction that relies on
these cannot be generalized easily to the quasi-setup.\footnote{However, also see \cite{HaNi99, Sch02}.}

\section{Candidate 2 : Weak Hopf Algebras}

Following the algebraists' intuition to avoid non-coassociativity and unnatural
ways to define duals, we move on to our second candidate: Weak Hopf algebras.
This time our objects are both algebras and coalgebras, associative and
coassociative, respectively, but the relationship between the two types of
structures is weakened. In particular we no longer force the coalgebra
structure to respect the unit of the algebra structure. In other words, we drop
the requirement that the comultiplication be unit preserving. We instead write:
\[ \Delta(1) = 1_{(1)} \otimes 1_{(2)} \neq 1 \otimes 1.\] 
\noindent
This, in turn, forces the counit to be \textit{weakly multiplicative}:
\[ \epsilon(xy) = \epsilon(x 1_{(1)}) \epsilon(1_{(2)} y). \]

Here is our precise definition:\footnote{We should point out that there is an alternative use in the literature for the term \emph{weak Hopf algebra}. Some authors use the term to refer to a bialgebra $B$ with a \emph{weak antipode} $S$ in ${\rm Hom}(B,B)$, i.e., ${id}*S*{id}={id}$ and $S*{id}*S=S$, where $*$ is the convolution product; see for instance \cite{Li98, LiDu02}. \cite{Wi02} briefly discusses how the two concepts are related.}
\begin{definition}
A (finite) \textbf{weak Hopf algebra} over $\kk$ is a finite-dimensional\footnote{We define only the finite-dimensional version here. In this context, it seems to be standard to restrict to finite dimensions because it is much easier to dualize the notion in that case. For a more detailed discussion of why finite-dimensionality is in general preferred, we refer the reader to \cite{BNS}.}
$\kk$-vector space $H$ with the structures of an associative algebra $(H,m,1)$, 
and a coassociative coalgebra $(H,\Delta,\epsilon)$ such that:
\begin{enumerate}
\item The comultiplication $\Delta$ is a (not necessarily unit-preserving)
homomorphism of algebras such that:
\[ (\Delta \otimes id) \Delta(1) = (\Delta(1) \otimes 1)(1 \otimes \Delta(1)) =
(1 \otimes \Delta(1)) (\Delta(1) \otimes 1);\]
\item The counit $\epsilon$ is a $\kk$-linear map satisfying the identity:
\[ \epsilon(fgh) = \epsilon(f g_{(1)})\epsilon(g_{(2)}h) =  \epsilon(f g_{(2)})\epsilon(g_{(1)}h) \]
\noindent
for all $f,g,h \in H$;
\item There is a linear map $S : H \rightarrow H$, called an \textbf{antipode},
such that for all $h \in H$, 
\begin{eqnarray*}
m(id \otimes S) \Delta(h) &=& (\epsilon \otimes id)(\Delta(1)(h\otimes 1));
\\
m(S \otimes id) \Delta(h) &=& (id \otimes \epsilon)((1 \otimes h)\Delta(1));
\\
m(m \otimes id)(S \otimes id \otimes S)(\Delta \otimes id)\Delta(h) &=&
S(h).
\end{eqnarray*}
\end{enumerate}
\end{definition}

We can easily see that the condition on the counit implies weak
multiplicativity. We also note that the defining equations for the antipode
given above would coincide with Equation \eqref{antipodedef} if $\Delta$ were
to preserve the unit. This would also imply that the counit would be a
homomorphism of algebras, and we would thus end up with a Hopf algebra. 

Weak Hopf algebras behave much better with respect to duality. In other words,
given a weak Hopf algebra $A$ over a field $\kk$, the space $Hom_{\kk}(A,\kk)$
can be given a natural weak Hopf algebra structure, using the canonical pairing
$\langle \quad , \quad \rangle\colon Hom_{\kk}(A,\kk) \times A\to \kk$.

Weak Hopf algebras in the above sense were introduced in \cite{BNS, BSz96, N, Sz96}. Later on it was observed that the paragroups of \cite{Oc88} and the face algebras of \cite{Hay93, Hay94, Hay98} fit nicely into this description. Many of the original constructions of weak Hopf algebras were motivated by applications to operator algebras, but many saw from early on their relationship to ``\textit{dynamical deformations of quantum groups.}"  We refer the reader to \cite{NV} for a more detailed exposition of weak Hopf algebras, and their relationship to various generalizations of the idea of quantum groups. This reference also contains discussion of various examples of weak Hopf algebras, one of which we present here.

Recall that group algebras were basic examples of Hopf algebras. We will look at
a nice generalization of these to find our weak Hopf algebras. In particular we
will need to generalize the notion of a \emph{group} via the
following:

\begin{definition}
A \textbf{groupoid over a set $X$} is a set $G$ together with the following
structure maps:
\begin{enumerate}
\item A pair of maps $s, t : G \rightarrow X$, respectively called the
\emph{source} and the 
\emph{target}.
\item A \emph{product} $m$, i.e. a partial function $m : G \times G \rightarrow
G$ satisfying the
following two properties:\footnote{The domain of $m$ is called the \emph{set of composable pairs}.}
\begin{enumerate}
\item $t(m(g,h)) = t(g)$, $s(m(g,h)) = s(h)$ whenever $m(g,h)$ is defined;
\item $m$ is associative: $m(m(g,h),k) = m(g,m(h,k))$ whenever the relevant
terms are defined. 
\end{enumerate}
\item An embedding $\epsilon : X \rightarrow G$ called the \emph{identity
section} such that
$m(\epsilon(t(g)),g) = g = m(g, \epsilon(s(g)))$ for all $g \in G$.
\item An \emph{inversion map} $\imath : G \rightarrow G$ such that
$m(\imath(g),g) =
\epsilon(s(g))$ and $m(g,\imath(g)) = \epsilon(t(g))$ for all $g \in G$.
\end{enumerate}
\end{definition}

Note that a groupoid over a singleton $X = \{x\}$ is a group.

We should point out here that there are other ways to define a groupoid. One way is to view it as a particular type of category $(X,G)$ with $X$ making up the set of objects, such that the morphisms (elements of $G$) are all invertible. We
choose not to consider this categorical description. However, simply by drawing some diagrams to represent the structure maps defined above, one can clearly see how the categorical notion may be derived easily.

Groupoids were first introduced in 1926, and since then found applications in
differential topology and geometry, algebraic geometry and algebraic topology,
and analysis. A very accessible introduction to groupoids can be found in
\cite{Wei2}. For more rigorous accounts one may refer to the bibliography
there. 

Now we start with a finite groupoid $G$ and consider the algebra $\kk G$. Here
we are considering the product $m(g,h)$ of two elements of $G$ to be defined as
in the groupoid itself (and in cases when it is not defined, we set the
relevant product equal to zero). We define the comultiplication $\Delta$, the
counit $\epsilon$ and the antipode $S$ on $G$ by:
\[ \Delta(g) = g \otimes g; \quad  \epsilon(g) = 1; \quad S(g) = \imath(g) \quad
\textmd{ for all } g \in G.\]
\noindent
Note that when $G$ is in fact a group, $\kk G$ is its group algebra and hence a
Hopf algebra.

Incidentally the dual of the groupoid algebra is again a weak Hopf algebra and
can be viewed as the function algebra on $G$. See \cite{NV} for more on this
and other examples.

Due to the inherent symmetry in their construction, weak Hopf algebras have
appealed to many researchers with strong algebraic preferences. Their theory
has been studied in detail and many algebraically natural constructions have
been generalized to their context, see for instance \cite{BNS, BSz00, Vec}. Even the curious fact that the trivial representation of a weak Hopf algebra may or may not be indecomposable (\cite[Prop.2.15]{BNS}) has interesting implications; see \cite{BSz00}.

Weak Hopf algebras have nice representation categories \cite{BSz00, NTV}. In
particular, representation categories of semisimple finite weak Hopf algebras
are (multi-)fusion categories. A \emph{fusion category} is a semisimple rigid monoidal category with finitely many simple objects and finite-dimensional homomorphism spaces, such that the endomorphism space of the unit object is one-dimensional. Relaxing the condition that the endomorphism space of the unit object be one-dimensional gives us the definition of a \emph{multi-fusion category}.
Fusion categories have been of much interest in the past few years. We refer
the reader to \cite{ENO} for a presentation of recent results about them, and
\cite{CE} for another accessible exposition. The reader may also find
\cite{BaKi} useful for background. 

It can be shown that a multi-fusion category can be viewed as the representation category of a (non-unique) semisimple weak Hopf algebra \cite{Hay99, Sz00}. In fact it turns out that many natural category-theoretical constructions can be restated in the language of weak Hopf algebras. Thus in \cite{ENO}, we see many results on (multi-)fusion categories stated and proved in terms of weak Hopf algebras. The authors of \cite{ENO} also show that a certain class of fusion categories can also be realized as the representation categories of finite-dimensional semisimple quasi-Hopf algebras. Hence we can see that the notion of a fusion category allows us to build a framework in which both quasi-Hopf algebras and weak Hopf algebras can be understood.

At this point, we move on to our third candidate, the Hopf algebroid, even
though our second one is still quite promising. We provide the reader with only
one reason, and that is generality. In \cite{BSz04}, B\"{o}hm and Szlach\'{a}nyi
showed that weak Hopf algebras with bijective antipodes are in fact Hopf
algebroids. In \cite{Xu99}, Xu constructed particular Hopf algebroids which
fully encoded the information of certain quasi-Hopf algebras associated with
solutions of the quantum dynamical Yang-Baxter equation. 
In other words, our third candidate will in fact provide us with a fine
structure which incorporates most of the interesting weak Hopf algebras and all
of the dynamical quantum groups realized as quasi-Hopf algebras, and then
offers us some more.

\section{Candidate 3 : Hopf algebroids}
 \label{HopfoidSection}
 
A natural question for an algebraist would be: What if we do not restrict
ourselves to commutative rings $\kk$? Can one develop the theory of bialgebras
in this more general setting? In 1977, Takeuchi \cite{Tak} described and
studied a new algebraic structure generalizing bialgebras to the noncommutative
setting. His structures were originally called \emph{$\times_A$-bialgebras}.
Then in 1996, when developing a geometrically motivated generalization of the
theory of quantum groups, Lu defined similar structures \cite{Lu1} and called
them \emph{bialgebroids}, to emphasize that the way these structures generalized bialgebras resembled the way Poisson groupoids generalized Poisson groups. A bit later, in 2001 \cite{Xu01}, Xu gave a similar definition and also implied the equivalence of these three notions. \cite{BM} provides a complete algebraic proof of this result. \cite{Bo1, BSz04, BM} all include a historical account of
these developments. In this note we will mainly be following \cite{Bo1} and
\cite{BSz04} for our definitions.

In short, a bialgebroid should be the natural extension of the notion of a
bialgebra to the world of groupoids. This then implies that a bialgebroid is no
longer an algebra, but a bimodule over a non-commutative ring. Here we will describe two symmetrical notions, that of a \emph{left bialgebroid} and a \emph{right bialgebroid}, first developed in \cite{KaSz}. We begin first with a \emph{left bialgebroid} $\mathcal{B}_L$, given by the following data:
\begin{enumerate}
\item Two associative unital rings: the \emph{total ring} $B$ and the
\emph{base ring} $L$.
\item Two ring homomorphisms: the \emph{source} ${s_L : L \rightarrow B}$ and
the \emph{target} ${t_L : L^{op} \rightarrow B}$ such that the images of $L$
in $B$ commute, making $B$ an $L-L$ bimodule denoted by $_LB_L$.
\item Two maps ${\gamma_L : B \rightarrow B_L \otimes {_LB}}$ and ${\pi_L : B
\rightarrow L}$, which make ${(_LB_L, \gamma_L, \pi_L)}$ a comonoid in
the category of $L-L$ bimodules. 
\end{enumerate}

We also need the following to be satisfied:
\begin{eqnarray*}
b_{(1)}t_L(l)\otimes b_{(2)} &=& b_{(1)}\otimes b_{(2)}s_L(l),
\\
\gamma_L(1_B) = 1_B \otimes 1_B
& and &
\gamma_L(b b') = \gamma_L(b) \gamma_L(b'),
\\
\pi_L(1_B) &=& 1_L,
\\
\pi_L(b s_L \circ \pi_L(b')) &=& \pi_L(b b') = \pi_L(b t_L \circ \pi_L(b'))
\end{eqnarray*}
\noindent
Note that in the above, we have made use of the Sweedler notation for ${\gamma_L : B \rightarrow B_L \otimes {_LB}}$: 
\[ \gamma_L(b) = b_{(1)} \otimes b_{(2)}. \]
\noindent
We also note that the source and target maps $s_L$ and $t_L$ may be used to define four commuting
actions of $L$ on $B$; these in turn give us in an obvious way the new
bimodules $^LB^L$, $^LB_L$, and $_LB^L$.

Similarly we can define a \emph{right bialgebroid} $\mathcal{B}_R$ using the
following data:
\begin{enumerate}
\item Two associative unital rings: the \emph{total ring} $B$ and the
\emph{base ring} $R$.
\item Two ring homomorphisms: the \emph{source} ${s_R : R \rightarrow B}$ and
the \emph{target} ${t_R : R^{op} \rightarrow B}$ such that the images of $R$
in $B$ commute, making $B$ an $R-R$ bimodule denoted by $^RB^R$. 
\item Two maps ${\gamma_R : B \rightarrow B^R \otimes {^RB}}$ and ${\pi_R : B
\rightarrow R}$, which make ${(^RB^R, \gamma_R, \pi_R)}$ a comonoid in
the category of $R-R$ bimodules. 
\end{enumerate}

As in the case of left bialgebroids, we need some conditions on $s_R, t_R, \gamma_R, \pi_R$ analogous to the conditions on $s_L, t_L, \gamma_L, \pi_L$:
\begin{eqnarray*}
s_R(r)b^{(1)}\otimes b^{(2)} &=& b^{(1)}\otimes t_R(r) b^{(2)},
\\
\gamma_R(1_B) = 1_B \otimes 1_B
& and &
\gamma_R(b b') = \gamma_R(b) \gamma_R(b'),
\\
\pi_R(1_B) &=& 1_R,
\\
\pi_R(s_R \circ \pi_R(b)b') &=& \pi_R(b b') = \pi_R(t_R \circ \pi_R(b)b')
\end{eqnarray*}
\noindent
Note that in the above, we have made use of a modified version of the Sweedler notation for ${\gamma_R : B \rightarrow B^R \otimes {^RB}}$: 
\[ \gamma_R(b) = b^{(1)} \otimes b^{(2)}. \]
\noindent
Also we can define three other bimodule
structures on $B$ using the source and the target, and denote them by $_RB_R$,
$^RB_R$, and $_RB^R$. 

These bimodule structures and the two notions of
bialgebroids are related as expected. For instance if $\mathcal{B}_L = (B, L,
s_L, t_L, \gamma_L, \pi_L)$ is a left bialgebroid, then its co-opposite is
again a left bialgebroid: $(\mathcal{B}_L)_{cop} = (B, L^{op}, t_L, s_L,
\gamma_L^{op}, \pi_L)$, where $\gamma_L^{op}$ is defined as $T \circ
\gamma_L$. (Here, as before, $T$ is the usual twist map, sending $a
\otimes b$ to $b \otimes a$). The opposite $(\mathcal{B}_L)^{op}$ defined by the data  $(B^{op}, L, t_L, s_L, \gamma_L, \pi_L)$ is a right bialgebroid. 

We refer the reader to \cite{BM} for several concrete examples of bialgebroids.
In this reference, Brzezinski and Militaru show a way to associate a
bialgebroid to a braided commutative algebra in the category of Yetter-Drinfeld
modules. Moreover they show that the smash product of a Hopf algebra with an
algebra in the Yetter-Drinfeld category is a bialgebroid if and only if the
algebra is braided commutative. Another interesting construction from \cite{BM}
gives a generic method of obtaining bialgebroids from solutions of the quantum
Yang-Baxter equation.

As mentioned above, Takeuchi's $\times_A$-bialgebras \cite{Tak}, Lu's
bialgebroids \cite{Lu1} and Xu's bialgebroids \cite{Xu01} were all shown to be
equivalent \cite{BM}, and all these are compatible with the definition above
which we took from \cite{Bo1, BSz04}. Currently, therefore, there is a universal
consensus on what should be accepted as the correct structure generalizing
bialgebras to the noncommutative base ring case. How, then, does one develop
the theory of Hopf algebras in this same setting? In the following paragraphs,
we briefly look at how several algebraists approached this problem. But first
we would like to answer a different question: Why were algebraists interested
in generalizing to the noncommutative base ring case in the first place?

In the mid 1990s, certain finite index depth 2 ring extensions from the theory
of Von Neumann algebras were shown to be related to some Hopf algebras
\cite{Szy}. A search started to find connections with more and more general
extensions. In 2003, some connections were discovered with the bialgebroids of
Takeuchi \cite{KaSz}. Thus the case which would correspond to some bialgebroid
with a generalized antipode in the noncommutative base ring case was naturally
of some interest. 

However, this problem proved to be somewhat complicated. The original antipode
suggested in \cite{Lu1} was not universally accepted, 
and various other formulations followed. We refer the reader to \cite{Bo1} for a
comparative study of these various antipodes.\footnote{A more category-theoretic discussion may be found in \cite{DaSt}.} The definition used in this note
is from \cite{BSz04}.

In
particular, to define a Hopf algebroid, we need two associative unital
rings $H$ and $L$, and set $R = L^{op}$. We consider a left bialgebroid
structure $\mathcal{H}_L = (H, L, s_L, t_L, \gamma_L, \pi_L)$ and a right
bialgebroid structure $\mathcal{H}_R = (H, R, s_R, t_R, \gamma_R, \pi_R)$
associated to this pair of rings. We require that $s_L(L) = t_R(R)$ and
$t_L(L) = s_R(R)$ as subrings of $H$, and:
\begin{eqnarray*}
(\gamma_L \otimes id_H) \circ \gamma_R &=& (id_H \otimes \gamma_R) \circ
\gamma_L, \\
(\gamma_R \otimes id_H) \circ \gamma_L &=& (id_H \otimes \gamma_L) \circ
\gamma_R. 
\end{eqnarray*}
The last ingredient is the antipode; this will be a bijection $S : H
\rightarrow H$, which will satisfy:\footnote{We require that the antipode $S$ be a bijection for the sake of simplicity. For a study of Hopf algebroids with antipodes which are not necessarily bijective, see \cite{Bo2}.}
\begin{eqnarray*}
S(t_L(l)ht_L(l')) &=& s_L(l')S(h)s_L(l), \\
S(t_R(r')ht_R(r)) &=& s_R(r)S(h)s_R(r')
\end{eqnarray*}
\noindent
for all $l,l' \in L$, $r,r' \in R$, and $h \in H$. In other words, we require
$S$ to be a twisted isomorphism simultaneously of bimodules $^LH_L \rightarrow
{_LH^L}$ and of bimodules $^RH_R \rightarrow {_RH^R}$. 

Our final constraint on the antipode $S$ is as follows: 
\begin{eqnarray*}
S(h_{(1)})h_{(2)} &=& s_R \circ \pi_R(h), \\
h^{(1)}S(h^{(2)}) &=& s_L \circ \pi_L(h)
\end{eqnarray*}
\noindent
for any $h \in H$. Once again, the subscripts and the superscripts on $h$ come from a
generalized version of the Sweedler notation which we use to define the
two maps $\gamma_L$ and $\gamma_R$: 
\begin{eqnarray*}
\gamma_L(h) &= h_{(1)} \otimes h_{(2)} &\in H_L \otimes {_LH} \\
\gamma_R(h) &= h^{(1)} \otimes h^{(2)} &\in H^R \otimes {^RH}.
\end{eqnarray*}

In this setup, then, we say that the triple $ \mathcal{H} =
(\mathcal{H}_L, \mathcal{H}_R, S)$ is a \emph{Hopf algebroid}.\footnotemark\,
\footnotetext{We should note here that some of the above information in our
definition is redundant. In fact one can start with a left bialgebroid
$\mathcal{H}_L = (H, L, s_L, t_L, \gamma_L, \pi_L)$ and an anti-isomorphism
$S$ of the total ring $H$ satisfying certain conditions, and from here can
reconstruct a right bialgebroid $\mathcal{H}_R$ using the same total ring such
that the triple $(\mathcal{H}_L, \mathcal{H}_R, S)$ is a Hopf algebroid. See
\cite{Bo1} for more details. An even more streamlined set of defining axioms is proposed in \cite{BoBr}.} 

It is easy to see how this symmetric definition, in terms of two bialgebroids
and a bijective map called an antipode, is analogous to the definition of a
Hopf algebra from a bialgebra and a bijective map called an antipode.
However, it may not be nearly as easy to see why the particular conditions
above describe the ``\emph{right antipode.}" We accept the definition given
above without further analysis; a sufficient discussion of this issue is
provided already in the introduction to \cite{Bo1}.

There are many examples of Hopf algebroids in recent literature.
Here we will describe an interesting example from \cite{BSz04}.

Let $\kk$ be a field of characteristic different from 2. Consider the group
bialgebra $\kk \Z_2$ presented as a left bialgebroid with the relevant
operations on the single generator $t$ as follows:
\[ t^2 = 1; \quad \Delta(t) = t \otimes t; \quad \epsilon(t) = 1. \]
\noindent
The source and the target maps are the natural ones: ${s, t  : \kk \rightarrow
\kk \Z_2}$ defined as $s(\lambda) = \lambda 1= t(\lambda)$.
Now if we define an antipode $S : \kk \Z_2 \rightarrow \kk \Z_2$ by setting
$S(t) = -t$, we obtain a Hopf algebroid. (Note that the given antipode is not a
Hopf algebra antipode).

The notions of integrals and duals for a Hopf algebroid have already been
studied \cite{BSz04, Bo2}. It turns out that the natural self-duality of weak
Hopf algebra structures is not available in this theory. In order to define
duals one needs a foray into the theory of integrals a la \cite{BSz04}. However
the representation categories are still quite nice. In \cite{Bo2}, it is shown
that the category of (left-left) Hopf modules of a Hopf algebroid $\mathcal{H} =
(\mathcal{H}_L, \mathcal{H}_R, S)$ is equivalent to the category of (right)
modules over the left base ring $L$.

Hopf algebroids can be used to describe the previous constructions discussed in
this note. In other words, as we briefly mentioned at the end of the previous
section, both quasi-Hopf objects and weak Hopf objects fit into the Hopf
algebroid picture with some modifications. In particular, it can be shown that
weak Hopf algebras with bijective antipodes are in fact Hopf algebroids
\cite{BSz04}. Moreover, we already know that using the framework of fusion
categories, we can describe quasi-Hopf objects in terms of weak Hopf algebras.
We can see a connection between quasi-Hopf objects and Hopf algebroids even
more directly, if we look at \cite{Xu99}. There, Xu constructs Hopf algebroids
which fully encode the information of certain quasi-Hopf algebras associated
with solutions of the quantum dynamical Yang-Baxter equation. Hence, it is
clear that there are some very interesting connections between these three
generalizations of the notion of Hopf algebras. 

\section{Short Interlude - Is this all that there is?}
\label{InterludeSection2}

In Section \ref{InterludeSection1}, we motivated our interest in the search for \emph{the correct generalization} of Hopf algebras by emphasizing the need for a Hopf-like object that can be used to develop sufficiently the theory of dynamical quantum groups. Indeed, several researchers already have investigated each of the three candidate structures we described so far with a view toward the theory of dynamical quantum groups. 

We discussed, already, how these three fit into the theory of dynamical quantum groups and also have seen how they relate to one another. 
In particular, the theory of (multi-)fusion categories seems to provide a fresh point of view which can become the right framework for understanding all of these structures.

Nonetheless, the search for other generalizations of Hopf algebras still goes on. In the rest of this note, we will introduce two more recently developed structures, also generalizing Hopf algebras. These have come up in contexts which are not immediately connected to the theory of dynamical quantum groups. However they both are structures which raise other new questions, and so they may be of interest to readers who are looking for general themes or for other new structures that may be useful for their own purposes. 

\section{A Category-Theoretical Approach: Hopf Monads}
\label{HopfMonadSection}

As mentioned in Section \ref{HopfAlgebraDefinitionsSection}, any finitely
generated commutative Hopf algebra over a field $\kk$ with characteristic 0 is
the function algebra of an affine group $G$ over $\kk$. Similarly a
cocommutative Hopf algebra generated by its primitive elements is the universal
enveloping algebra of a Lie algebra $\g$. Readers may find material on these
well-known results in classics like \cite{Ho81, MiMo65} and in more modern
texts like \cite{EtScbook, FeRi05}.

Both these results are in fact particular instances of \emph{Tannaka-Krein
duality}. Traditionally, the origins of Tannaka-Krein theory are attributed to
Groethendieck. We refer
the reader to \cite{HeRo70} for a comprehensive account relating the original
works of Tannaka and Krein with more modern treatments. The first modern references in the subject that make extensive use of category theory are \cite{SR, DM}. Expositions of some
Tannaka-Krein type theorems presented in the flavor closest to the perspective
of this note may be found in \cite{EtScbook}. For more details on the modern
approach, with an emphasis on Hopf algebras and monoidal categories, the reader may find
\cite{Maj95} useful.

In the language of categories, the philosophy underlying the various results
that can be gathered under the heading of Tannaka-Krein theory can be stated as
follows: it should be possible to view any \emph{nice} category $\mathcal{C}$
as the representation category of \emph{some algebraic structure}. This vague
assertion becomes accurate mathematics when one chooses appropriate
descriptions for the italicized phrases.

For instance, one precise formulation of
the statement of Section \ref{HopfAlgebraDefinitionsSection} about commutative
Hopf algebras is as follows:

\begin{theorem} \cite{DM, EtScbook}
Let $\mathcal{C}$ be a symmetric tensor category defined over an algebraically
closed field $\kk$ of characteristic zero. Suppose that $\mathcal{C}$ is
abelian and has finitely many indecomposable objects. Suppose, in addition,
that there is an exact faithful tensor functor $\mathbf{F} : \mathcal{C}
\rightarrow Vec_{\kk}$. Then the group $G$ of tensor automorphisms of
$\mathbf{F}$ is a finite group and $\mathcal{C} \simeq Rep(G)$.
\end{theorem}
In literature, one can find similar results for more general Hopf algebras
\cite{Sch, Ulb} and quasi-Hopf algebras \cite{Maj92}. In fact, the search for
ever more general Hopf-like objects always includes category-theoretic
investigations which look at the representation categories of these objects and
attempt to describe them as some nice types of categories.

Our fourth candidate for the correct generalization of Hopf algebras fits
perfectly into the spirit of this categorical approach to Hopf algebras. This is the so-called \emph{Hopf monad} structure.\footnote{We note here that there is at least one other category-theoretic Hopf-like object studied in the recent years. We do not focus on the relevant constructions here; we simply refer the reader to \cite{DaMCSt, DaSt} where the ideas are developed in great detail. Some connections to the picture from Section \ref{HopfoidSection} may be found in \cite{BSz04, Bo3}.}
Very briefly a \textbf{Hopf monad} is a Hopf-like object in a general category.

More precisely, we begin with a monoidal category $\mathcal{C}$, and assume that
it is \emph{rigid}; in other words, we assume that every object in
$\mathcal{C}$ has a left dual and a right dual. We consider first a
\emph{monad} in the sense of \cite{MacBook}. In short a \textbf{monad} on
$\mathcal{C}$ is an algebra in the monoidal category $End(\mathcal{C})$. More
specifically a monad is an endofuctor $\mathbf{T}$ of $\mathcal{C}$ equipped with two functorial morphisms $\mu : \mathbf{T}^2 \rightarrow \mathbf{T}$ and 
$\eta : 1_{\mathcal{C}} \rightarrow \mathcal{C}$ which satisfy certain conditions analogous to those describing the product and the unit of an algebra. 

The appropriate definition in this context for the notion generalizing
bialgebras was introduced by Moerdijk in \cite{Moe02} and asserts that a \textbf{bimonad} is a monad $\mathbf{T}$ which is also comonoidal. In the precise statement of the definition, there are functorial morphisms playing the analogous roles of the coproduct and the counit. The notion of a
bimonad is not self-dual, but defining duals is not too complicated.

To make a bimonad $\mathbf{T}$ into a Hopf monad, we only require the introduction of two new functorial morphisms, the so-called \emph{left antipode} and the \emph{right antipode}, which make use of the left and right duals of the
category $\mathcal{C}$ to encode the left and right duals in the category of
$\mathbf{T}$ modules. Thus a \textbf{Hopf monad} is a bimonad with two functorial morphisms, a left antipode and a right antipode.

Hopf monads in the above sense\footnote{We need to remark on our specific terminology here. The term ``\emph{Hopf monad}" appeared earlier, in \cite{Moe02}. In particular, the ``\emph{Hopf monads}" of \cite{Moe02} are precisely the structures that we here called ``\emph{bimonads}." We are following \cite{BrVi, Sz03} with this choice of terminology, because it seems quite natural linguistically to reserve the Hopf term for bimonads with antipode. However, for readers familiar with operad theory, the terminology choice of \cite{Moe02} may seem more reasonable; see also \cite{Moe01}.} 
were introduced and studied in detail in \cite{BrVi}. The authors' motivation there stems from the various topological invariants constructed via
methods of what is now called quantum topology \cite{Kau93}. The original
construction by Reshetikhin and Turaev \cite{ReTu} 
of invariants for $3$-manifolds using quantum $sl_2(\C)$ and related (ribbon)
Hopf algebras has been generalized in recent years to constructions in more and
more general settings. The authors of \cite{BrVi} study Hopf monads as a
generalization of Hopf algebras to categories with no braiding so as to provide
the ultimate framework to understand all these newer invariants in a uniform
manner. 

Any finite-dimensional Hopf algebra $H$ over $\kk$ can be used to construct a
Hopf monad on the category of finite-dimensional $\kk$-vector spaces. More
generally, Hopf monads can be constructed using any rigid monoidal category
with an underlying algebraic structure. Thus, the Hopf monad concept brings
together all of the earlier structures we studied in this note: quasi-Hopf
algebras, weak Hopf algebras and Hopf algebroids all can be used to construct
Hopf monads in their relevant representation categories.\footnote{For a clarification of the relationship between bimonads and bialgebroids, see \cite{Sz03}. The quasi- and weak Hopf algebra cases are explicitly discussed in \cite{BrVi}.}

\section{A (Poisson-)Geometric Approach: Hopfish Algebras}
\label{HopfishSection}

Finally we briefly describe a fifth Hopf-like object which generalizes the
notion of a Hopf algebra.
This object, playfully called a \emph{hopfish algebra} by its creators, was
introduced first in \cite{TaWeZh} and applied to the problem of describing
modules of irrational rotation algebras in \cite{BlTaWe}. The ideas in its
development fit in the context of (Poisson) geometry; the reader may find
interesting background and motivational discussions in the recent survey
\cite{BlWe} of algebraic structures in Poisson geometry. 

Hopfish algebras live in a category $\mathcal{A}$ whose objects are $\kk$-unital
algebras and whose morphisms are bimodules. Then as the authors of
\cite{TaWeZh} show, the notions of comultiplication, counit, and the antipode
can all be developed in this framework, as particular types of bimodules (i.e.
morphisms in the new category). 

More precisely, we start with the notion of a \emph{sesquiunital sesquialgebra},
which is the relevant generalization of bialgebras in this context:  A
\textbf{sesquiunital sesquialgebra} over a commutative ring $\kk$ is a unital
$\kk$-algebra $B$ equipped with a $(B \otimes B, B)$-bimodule ${\bf \Delta}$
(the \textbf{coproduct}), and a $(\kk,B)$-bimodule (ie a right $B$-module)
${\bf \epsilon}$ (the \textbf{counit}), with the following properties:
\begin{enumerate}
\item (\emph{coassociativity axiom}) The two $(B \otimes B \otimes B, B)$ bimodules
${(B \otimes {\bf \Delta})\otimes_{B \otimes B} {\bf \Delta}}$ and $({\bf \Delta} \otimes B)\otimes_{B \otimes B} {\bf \Delta}$ are isomorphic; and
\item (\emph{counit axiom}) The two $(B, B)$ bimodules
$({\bf \epsilon} \otimes B) \otimes_{B \otimes B} {\bf \Delta}$ and $(B \otimes
{\bf \epsilon}) \otimes_{B \otimes B} {\bf \Delta}$ are both isomorphic to $B$.
\end{enumerate}

To move on to the Hopf-like objects one needs an antipode appropriate for the
context. In order to define the right structure that should correspond to the
antipode, we begin with a \emph{preantipode}. In particular a
\textbf{preantipode} for a sesquiunital sesquialgebra $B$ over $\kk$ is a left
$(B \otimes B)$-module ${\bf S}$ together with an isomorphism of its $\kk$-dual
with the right $(B \otimes B)$-module $Hom_B({\bf \epsilon},{\bf \Delta})$. One
can view ${\bf S}$ as a $(B, B^{op})$-bimodule, and hence an
$\mathcal{A}$-morphism in $Hom(B, B^{op})$. If a preantipode ${\bf S}$ on $B$
is a free left $B$-module of rank $1$ when considered as a $(B,
B^{op})$-bimodule, then it is called an \textbf{antipode} for $B$. A
sesquiunital sesquialgebra $H$ equipped with an antipode ${\bf S}$ is called a
\textbf{hopfish algebra}. 

To see the standard examples, we need to begin with a well-known functor of
Morita theory. This functor, called \textbf{modulation} in \cite{TaWeZh}, goes
from the category of $\kk$-unital algebras with algebra homomorphisms as
morphisms, to the category $\mathcal{A}$. Using this functor, one can see that
the modulation of biunital bialgebras will be basic examples of  sesquiunital
sesquialgebras and the modulation of Hopf algebras will be basic examples of
hopfish algebras. In \cite{TaWeZh} the authors also show that quasi-Hopf and
weak Hopf algebras algebras can be modulated to yield hopfish algebras.

For more details on hopfish algebras we refer the reader to \cite{TaWeZh, BlTaWe, BlWe}. Readers who prefer the more general setting of higher category theory may find that \cite{DaMCSt, DaSt} can perhaps complement the above references on hopfish algebras. 

\section{Conclusion - Where do we go next?}
\label{ConclusionSection}

In this note, we focused on five structures that have been introduced in the recent
years as generalizations of Hopf algebras. 
One thing is clear:
The notion of Hopf algebras proved useful in so many diverse ways that mathematicians of
many different persuasions decided that the search for \emph{the correct
generalization} was quite an important task. More significantly, a lot of new
and interesting mathematics came up during these investigations. 

Some readers may also have noticed along the way a few open problems that are looking for solutions. For instance, the question whether the Hopf algebroids from \cite{Lu1} 
are a proper subclass of the Hopf algebroids of Section \ref{HopfoidSection} or not is still not settled. The example we gave at the end of Section \ref{HopfoidSection} would not be a Hopf algebroid in the sense of \cite{Lu1}. So far, we do not have an example which would be a Hopf algebroid for \cite{Lu1} and would not be a Hopf algebroid in our sense, but it is not certain that we will never find one. This is an interesting question, as the Hopf algebroids in \cite{Lu1} are quite natural structures themselves. 

There are also many open questions about (multi-)fusion categories and hence about (weak) Hopf algebras. One of the most fascinating questions involves the existence (or non-existence) of finite dimensional semisimple Hopf algebras with non-group theoretical representation categories. For details on this and many other open questions about (multi-)fusion categories we refer the reader once again to \cite{ENO}.

The relationship between representation categories of all the various generalizations of Hopf algebras also needs to be studied further. Indeed, a lot of the pieces are out there, but there is still more work to be done. 

In this note, we pointed out some such connections. We also provided some details on the motivations of the mathematicians who developed the relevant structures. Along the way, we gave our own reasons for being interested in generalizations of Hopf algebras. We hope that at this point, the readers have already decided if any of the five structures fits their own particular needs or mathematical
inclinations. If not, they at least have an idea of what leads to follow to
decide on this matter for themselves.

\end{document}